\documentclass[lettersize,journal]{IEEEtran}
\usepackage{amsmath,amsfonts}
\usepackage{algorithm}
\usepackage{array}
\usepackage[caption=false,font=scriptsize,labelfont=sf,textfont=sf]{subfig}
\usepackage{textcomp}
\usepackage{stfloats}
\usepackage{url}
\usepackage{verbatim}
\usepackage{graphicx}
\usepackage{booktabs}
\usepackage{cite}
\usepackage{amsthm}

\usepackage{multirow}

\usepackage{algpseudocode}

\newtheorem{Prp}{Proposition}

\newtheorem{Def}{Definition}

\begin{document}

\title{Transient Stability-Constrained OPF: Neural Network Surrogate Models and Pricing Stability}

\author{M. Garcia\IEEEauthorrefmark{1}, N. LoGiudice\IEEEauthorrefmark{2}, R. Parker\IEEEauthorrefmark{1}, R. Bent\IEEEauthorrefmark{1}

\thanks{\IEEEauthorrefmark{1} Los Alamos National Laboratory}
\thanks{\IEEEauthorrefmark{2} Electrical and Computer Engineering, Texas A\&M University}
}



\maketitle

\begin{abstract}
A Transient Stability-Constrained Optimal Power Flow (TSC-OPF) problem is proposed that enforces frequency stability constraints using Neural Network (NN) surrogate models.  NNs are trained using a novel model-driven active sampling algorithm that iteratively generates NN training data located near the stability boundary and contained within the feasible set of the Alternating Current Optimal Power Flow (AC-OPF) problem.  In the context of wholesale electricity markets, pricing structures are analyzed along with their dependencies on the selected input features to the NN surrogate model.  An important insight identifies a trade-off between the accuracy of the NN surrogate model and sensible locational pricing structures.  NN surrogate models for frequency stability are validated by ensuring the resulting TSC-OPF solution is stable over randomly generated load samples using a small Hawaii test case.  The proposed TSC-OPF problem is shown to significantly enhance frequency stability at low computational cost and low financial cost to the system.  For certain selections of NN inputs, the TSC-OPF problem is able to stabilize all load scenarios for which the solution to the AC-OPF problem resulted in instability.  
\end{abstract}


\section{Introduction}

Modern power systems have long relied on Economic Dispatch (ED) optimization models to assist in making decisions about which generators to utilize on hourly time scales. These optimization models assume steady-state conditions and their solutions have historically not required significant adjustment when deployed and operated in practice.  This has started to change in recent years, in part due to increasing penetrations of Inverter-Based Resources (IBRs) and renewable energy sources. Typically, these resources do not respond to frequency, have replaced resources that do respond to frequency, and have fostered environments that are prone to issues with frequency dynamic stability.  This requires ED problems to model adjustments on faster timescales and to account for system transients for ensuring that the dispatch solutions are reliable under such conditions. 

To address this problem, some Independent System Operators (ISOs) have introduced novel reserve products into their electricity markets to accommodate droop control and improve frequency dynamics~\cite{mago2020new,operator2017fast,meng2019fast}.  In this context, constraints on reserves are used to model stability requirements in ED problems, e.g.~\cite{li2018design,liu2018participation}. While these methods are an excellent start, fully addressing stability in ED generally remains an open challenge for a few reasons.  First, accurate power system dynamics models are represented by intricately detailed, non-linear Differential Algebraic Equations (DAEs) and appropriate stability criteria are difficult to derive.  Second, incorporating stability constraints into ED models makes these optimization problems more difficult to solve.  Finally, existing ED formulations are tightly coupled to modern energy markets and incorporating stability constraints could have poorly understood effects on markets and pricing. 

We address this problem by developing a novel framework for introducing stability criteria into ED problems in the form of a demonstrably accurate NN binary classifier that doesn't significantly increase the computational effort required to solve the ED problem.   We additionally show how such a model naturally fits into existing pricing mechanisms. Finally, one of the major challenges in developing NN surrogates for stability is the need for large amounts of measurement data, supplemented with computationally demanding simulation data, for training. To overcome this challenge, a model-driven, active sampling method is developed that dramatically reduces the number of samples needed for a high-quality model. Together, these contributions yield a computationally tractable model for solving the Transient Stability-Constrained Optimal Power Flow (TSC-OPF). 

\subsection{Literature Review} 
The ED optimization problem is used to determine a dispatch schedule for generators minutes to hours in advance and is also used to set electricity prices and clear the wholesale electricity market. ED problems typically include reliability constraints in the form of N-1 criteria that ensure there exists a feasible post-contingency dispatch after the loss of any single element in the system.  These constraints typically use quasi-steady state models to represent the transition from the pre-contingency dispatch to the post-contingency generator outputs, e.g. droop and voltage response~\cite{aravena2023recent}.  However, in practice the power system undergoes a dynamic (transient) response after a contingency that could result in instabilities or other issues not captured by quasi-steady state models; thereby creating the need for models of transient phenomena.

Embedding power system stability constraints into ED problems is studied within the literature.  There are generally three approaches to formulating these constraints in ED problems.  The first approach discretizes the DAEs describing the power system dynamics and directly embeds the resulting algebraic constraints into the ED problem~\cite{gan2000stability,jiang2010enhanced}.  This approach is computationally burdensome and is only tractable for small systems with simplified dynamic models. The second approach derives stability constraints from first principles~\cite{chavez2014governor,doherty2005frequency}.  The resulting stability constrained optimization problems are tractable to solve, but require many simplifying assumptions that are not appropriate for realistic power system applications.  The final approach runs detailed dynamic simulations to generate data used to train surrogate models that represent stability.  These stability surrogate models are then embedded into ED problems~\cite{ahmadi2013two,gutierrez2010neural}.  This approach has received some recent attention and is becoming more attractive as machine learning tools and capabilities have undergone significant advancements. This paper's approach falls in this last category.

Focusing on stability constraints derived from first principles, reference \cite{zhou2024review} provides a recent review of frequency stability constraints for ED problems.  Recent advancements in this area introduce details into stability constraints such as regional frequency deviations or novel ancillary services for improving frequency response~\cite{badesa2020conditions1,badesa2020conditions2,garcia2019real,trovato2018unit,zhang2020modeling}.  These approaches are considered computationally tractable and are often well studied in the context of ISO operations, electricity markets, and pricing~\cite{badesa2020pricing,badesa2022assigning,garcia2022Hawaii,garcia2021requirements}.  However, these works make significant assumptions including the use of over simplified turbine governor models, often assuming that droop control follows a predetermined piecewise linear trajectory, and simplified power flow models that neglect reactive power. Furthermore, this line of work neglects detailed dynamic models that are considered standard in power system simulation software, e.g. dynamic machine and exciter models specific to different generator types~\cite{PowerWorld,kundur1994power}. Despite the reliance on strong potentially unrealistic assumptions, these works typically do not validate ED solutions with high-fidelity dynamic simulations and it is unclear if these solutions ensure stability in practice.

In contrast, relatively few works have incorporated data-driven surrogate models for frequency stability into ED problems.  The most recent work has embedded NN surrogate models into Unit Commitment (UC) problems formulated as a Mixed-Integer Linear Program (MILP).  These problems typically use Rectified Linear Unit (ReLU) activation functions to fit into a MILP framework.  Reference \cite{zhang2021encoding} uses an active sampling algorithm to train a NN surrogate constraint that ensures the frequency will not fall below a specified lower frequency bound.  Reference \cite{wu2023transient} uses an input-convex NN to represent a stability index defined by the maximum rotor angle difference between generators. Another approach uses NN-based stability surrogate constraints in Alternating Current Optimal Power Flow (AC-OPF) problems represented as Non-Linear Programs (NLPs), i.e., TSC-OPF.  Reference \cite{su2021deep} uses a deep belief network to represent a transient stability function and solves the TSC-OPF problem using a genetic algorithm.  References \cite{ahmadi2013two} and \cite{gutierrez2010neural} use NN surrogates with smooth activation functions to represent a transient stability boundary and they solve the TSC-OPF problem using gradient-based optimization algorithms.

In implementation settings, the ISO in Texas has recently implemented stability constraints in the form of reserve requirements.  Assuming fixed inertia levels, a linear structure is imposed for these reserve requirements and the parameters of the resulting linear constraints are determined empirically using detailed dynamic simulations~\cite{li2018design,liu2018participation}.  This approach is remarkably similar to training NN surrogate models, but with a simpler mathematical structure.  For this reason, it is anticipated that NN surrogate approaches hold promise and are ready for practical implementations.  However, existing NN surrogate models have not been demonstrated in such practical settings as they do not account for economic implications.  For example, ISOs will likely enforce stability constraints as reserve requirements, where the reserve of interest procures generator headroom for resources to provide ancillary services that improve primary frequency control, e.g. droop control. Most importantly, the input features of the NN surrogate models heavily impact pricing structures in the context of electricity markets, a topic that has not yet been addressed in the literature.

\subsection{Contributions}
This paper develops a novel active sampling algorithm that differs from that in~\cite{zhang2021encoding} and related work in two important ways.  First, the proposed algorithm generates training samples within a model-driven, iterative process (active sampling), which uses the optimization of the OPF to dynamically generate training data and update the NN.  In contrast, the active sampling algorithm developed in \cite{zhang2021encoding} generates \emph{all} samples apriori.  Their active sampling procedure then chooses which of these finite number of samples to use during NN training. Second, the optimization problem solved during each iteration of our proposed algorithm has a continuous feasible set representing the feasible set of the AC-OPF problem.  In contrast, the algorithm from \cite{zhang2021encoding} solves a discrete optimization problem that simply optimizes over the finite set of available samples that were generated apriori.  We found these two improvements to be necessary in our setting where stability depends significantly on continuous dispatch variables, see Section \ref{Sec:GovLimits} regarding why stability depends significantly on the continuous dispatch variables. 

The proposed TSC-OPF problem is validated over randomly generated load samples by checking stability of the optimal dispatch.  This illustrates robustness of the proposed NN surrogate models and effectiveness of the proposed active sampling algorithm.  Indeed, previous work on the topic of NN surrogate models for power system stability do not sufficiently validate the optimal solutions of the TSC-OPF or TSC-UC problems over different problem instances.  Although it is common in previous work to validate NN accuracy over randomly sampled NN inputs, those randomly sampled NN inputs do not represent optimal solutions of the TSC-OPF or TSC-UC problems.  Instead, it is important to ensure that the optimization problems are not identifying optimal solutions located in poorly trained areas of the NN input space.  This is an important differentiator between our work and other work on the topic of TSC-OPF \cite{ahmadi2013two,gutierrez2010neural}.  

Finally, this paper investigates different input features to the NN surrogate model and their impact on pricing structures.  In this context it is observed that the most accurate choice of inputs result in discriminatory pricing structures that aren't well accepted by current electricity market practices.  However, the inputs can be chosen in a way that results in both high accuracy and uniform pricing structures.

\section{Transient Stability-Constrained OPF}\label{SecII}
This section formulates the general TSC-OPF optimization problem and provides a reformulation of the stability constraints using an explicit NN representation.  The parameters of the stability constraint are left general for study of different choices of input features.

\subsection{General Problem Formulation}
The total number of buses in the system is $n$.  Let $\bold{p}\in\mathbb{R}^n$ and $\bold{q}\in\mathbb{R}^n$ represent vectors of real and reactive power injections at each bus.  The set of generators is denoted by $\mathcal{G}$ and the total number of generators in the system is $m$.  Let $\bold{g}_j\in\mathbb{R}$, $\bold{r}_{j}\in\mathbb{R}$, and $\bold{h}_{j}\in\mathbb{R}$ represent the real generation, reactive generation, and reserve provided by generator $j$.  The matrix $M\in\mathbb{R}^{n\times m}$ maps each generator to its corresponding bus, where $M_{ij}=1$ if generator $j$ is located at bus $i$.   The cost of generator $j$ is assumed to be convex and is denoted by $C_j(\!\bold{g}_{j}\!,\!\bold{r}_{j}\!)$.  The generator model is general enough to accommodate market participants that consume electricity; however, this model is referred to as a generator model.  Electricity demand (or load) is modeled as fixed and is aggregated at each bus where $\bold{d}\in\mathbb{R}^n$ and $\bold{l}\in\mathbb{R}^n$ represent the real and reactive power load at each bus in the system.

We present an abstract TSC-OPF problem that accommodates many different types of stability constraints, different models of the transmission network, and different generator models.  The TSC-OPF problem minimizes generator dispatch cost and is written as follows:
\begin{subequations}\label{TSCOPF}
\begin{align} 
& \min \ \ \  \underset{j \in \mathcal{G}}{\textstyle\Sigma} C_j(\!\bold{g}_{j}\!,\!\bold{r}_{j}\!)\tag{\ref{TSCOPF}}\\
&(\boldsymbol{\lambda}) \ \ \ \  \bold{d}-M\bold{g}+\bold{p}= \mathbf{0} & 
\label{rpowerbalance}\\
&(\boldsymbol{\mu}) \ \ \ \  \bold{l}-M\bold{r}+\bold{q}= \mathbf{0} & 
\label{qpowerbalance}\\
&(\boldsymbol{\gamma}) \ \ \ \  f(\bold{g},\bold{d},\bold{r},\bold{l},Z\bold{h},\bold{p})\geq c & \label{stabilityconstraint}\\
&\hspace{30pt}(\!\bold{g}_{j}\!,\bold{r}_{j}\!,\bold{h}_{j}\!)\in\mathcal{X}_j \hspace{5pt}\forall j \in \mathcal{G}\label{GenPrivateConstraint}\\
&\hspace{30pt}(\bold{p},\bold{q})\in\mathcal{T},\label{TransmissionConstraint}
\end{align}
\end{subequations}
where \eqref{rpowerbalance} represents real-power balance at each bus, \eqref{qpowerbalance} represents reactive-power balance at each bus, and \eqref{stabilityconstraint} represents a general stability constraint.   Generator constraints and transmission constraints are enforced implicitly in \eqref{GenPrivateConstraint} and \eqref{TransmissionConstraint} respectively.  The Lagrange multipliers appear to the left of each constraint in parenthesis.  

The NN surrogate model function is denoted by $f$ and maps a set of input features (variables and parameters) to a \mbox{$\rho$-dimensional} space.  The output of this function is restricted by the constant $c\in\mathbb{R}^{\rho}$.  In this paper, this function outputs a scalar, e.g. $\rho=1$.  The inputs of $f(\cdot)$ specified in \eqref{stabilityconstraint} represent a super-set of the inputs considered in this paper.  This leaves the TSC-OPF formulation general to the input features used to train the NN. Notice that the load parameters $\bold{l}$ and $\bold{d}$ are inputs to the function $f(\cdot)$ even though they are not optimization variables.  Using these inputs allows the NN to be more accurate for different TSC-OPF problem instances defined by different load samples.  Finally, zonal aggregation of reserve is denoted $Z\bold{h}$ where $Z\in\mathbb{R}^{z\times m}$ is a matrix of ones and zeros where each element $Z_{ij}$ is one if generator $j$ is contained within zone $i$.

The AC-OPF problem is defined as \eqref{TSCOPF} with constraint \eqref{stabilityconstraint} removed.

\subsubsection{Transmission Constraints}
In this formulation, the transmission constraints are enforced implicitly by the set $\mathcal{T}$.  
Our work uses the AC transmission constraints with voltages expressed in polar coordinates as outlined in \cite{cain2012history}.  In this context $\mathcal{T}$ is non-convex and is described by constraints that include only continuously differentiable functions.   

\subsubsection{Generator Model}
The generator constraints are enforced implicitly with the set $\mathcal{X}_j$, for generator $j$. This set is left general to accommodate other types of market participants.  
 In this paper, the generator constraints are represented as follows:
\begin{subequations}\label{GenSet}
\begin{align}
\hspace{-5pt}\mathcal{X}_j:=\left\{\right.& (\bold{g}_j,\bold{r}_j,\bold{h}_j)\in \mathbb{R}\!\times\!\mathbb{R}\!\times\!\mathbb{R}_+: \tag{\ref{GenSet}}\\
&\hspace{-21pt}\bold{\underline{g}}_j\leq \bold{g}_j \leq \bold{\bar{g}}_j \\
&\hspace{-21pt}\bold{\underline{r}}_j\leq \bold{r}_j \leq \bold{\bar{r}}_j \\
&\hspace{-21pt}\bold{h}_j = \min(\bold{\bar{h}}_j,\bold{\bar{g}}_j-\bold{g}_j)\left.\right\},\label{nonsmoothconstraint}\end{align} \end{subequations}
where $\bold{\underline{g}}_j$ and $\bold{\bar{g}}_j$ represent the minimum and maximum real-power output and $\bold{\underline{r}}_j$ and $\bold{\bar{r}}_j$ represent the minimum and maximum reactive-power output for generator $j$.  The constant $\bold{\bar{h}}_j$ models the maximum reserve a generator $j$ is capable of providing.  Equation \eqref{nonsmoothconstraint} represents a non-smooth equality constraint that can be removed if reserve is not included as an input to the stability function $f(\cdot)$, which is the case in many of our numerical results.

Generators are also modeled as profit maximizers.  Generator $j$ is paid for generating real-power at a price $\pi_j$, reactive-power at a price $\sigma_j$, and reserve at a price $\alpha_j$.  Each generator is assumed to be a price-taker and thus prices are assumed fixed from their perspective.  The profit maximization problem for generator $j$ is as follows:
\begin{subequations}\label{ProfMax}
\begin{align} 
& \max \ \ \  \pi_j \bold{g}_j + \sigma_j \bold{r}_j + \alpha_j \bold{h}_j - C_j(\bold{g}_{j}\!,\!\bold{r}_{j}\!)\tag{\ref{ProfMax}}\\
&\hspace{35pt}(\!\bold{g}_{j}\!,\bold{r}_{j}\!,\bold{h}_{j}\!)\in\mathcal{X}_j.
\end{align}
\end{subequations}

\subsection{Explicit NN Representation}
A NN can be expressed mathematically as an iterative application of affine transformations that are passed through element-wise activation functions.  Each successive application of activation functions represents a different layer in the NN.  This representation of a NN is expressed as follows:
\begin{subequations}
\begin{align}
&\bold{y}^{(k)}=W^{(k)}\bold{x}^{(k-1)}+\bold{b}^{(k)}  \hspace{5pt} \forall k\in[1,\kappa],\label{explicit2}\\
&\bold{x}^{(k)}=\sigma_k(\bold{y}^{(k)}) \hspace{5pt} \forall k\in[1,\kappa-1],\label{explicit3}
\end{align}\end{subequations}
where $\bold{x}^{(0)}$ represents the input to the NN, the layers are indexed by $k$, and the element-wise activation function used for layer $k$ is denoted $\sigma_k(\cdot)$.  There are $\kappa$ layers in the NN each having potentially different number of neurons.  The number of neurons in layer $k$ is denoted $\beta_k$ and the matrix $W^{(k)}\in\mathbb{R}^{\beta_k\times\beta_{k-1}}$ and vector $b^{(k)}\in\mathbb{R}^{\beta_k}$ are used to define the affine transformations.

This mathematical expression of a NN is directly implemented into the TSC-OPF problem formulation as algebraic constraints using MathOptAI.jl~\cite{MathOptAI2024}.  In this context the input of the NN represents the chosen inputs of the function $f(\cdot)$, e.g. $\bold{x}^{(0)}=[\bold{g},\bold{d},\bold{r},\bold{l},Z\bold{h},\bold{p}]$.  This explicit formulation replaces constraint \eqref{stabilityconstraint} with constraints \eqref{explicit2} and \eqref{explicit3} along with the following inequality constraint on the NN output:
\begin{equation}
\bold{y}^{(\kappa)}\geq c.\label{explicit1}
\end{equation}

\section{Choice of NN Input and Pricing Implications}
One of the challenges with introducing NN surrogate constraints into ED is that the literature has not yet addressed how these constraints impact prices and market clearing.  This section introduces a methodology for formally evaluating the pricing implications of \eqref{TSCOPF} on the basis of the input features of \eqref{stabilityconstraint}.  In short, the specific choice of inputs impacts both the accuracy of the NN surrogate model and the pricing structure in the context of an electricity market.  This section develops a pricing structure that ensures the generator dispatch will maximize profits of each individual generator.  These prices encourage a generator to follow their dispatch instructions, effectively aligning incentives of all market participants. Additionally, these prices are shown to be locational for some choices of inputs and for other choices of inputs these prices are shown to be specific to each individual generator.  Locational prices are referred to as being \emph{uniform} and prices specific to each generator are referred to as being \emph{discriminatory}. 

\subsection{Defining Prices}\label{Sec:Pricing}

First, assume the TSC-OPF problem can be solved to a solution that satisfies the Karush–Kuhn–Tucker (KKT) conditions.  
\Def{\normalfont \label{KKTdef}A \emph{KKT point} is denoted $({\bold{g}}^{\star}\!\!,{\bold{r}}^{\star}\!\!,{\bold{h}}^{\star}\!\!,{\bold{p}}^{\star}\!\!, {\bold{q}}^{\star}\!\!;{\lambda}^{\star}\!\!,{\mu}^{\star}\!\!,{\gamma}^{\star})$ and satisfies the KKT conditions for the TSC-OPF problem \eqref{TSCOPF}.}\normalfont

The following prices are proposed for generator $j$ at bus $i$ where bus $i$ is contained within zone $k$:
\begin{subequations}
\begin{align}
&\pi_j= \lambda_i^{\star}+\nabla_{\bold{g}_j}\!f^{\star} \gamma^{\star}\label{pprice}\\
&\sigma_j= \hat{\mu}_i^{\star}+\nabla_{\bold{r}_j}\!f^{\star} \gamma^{\star}\label{qprice}\\
&\alpha_j= \nabla_{\!Z_k \bold{h}} f^{\star} \gamma^{\star}\label{hprice}
\end{align}\end{subequations}
where $\nabla_{\bold{g}_j}\!f^{\star}$, $\nabla_{\bold{r}_j}\!f^{\star}$, and $\nabla_{\!Z_k \bold{h}} f^{\star}$ represent the gradient of $f(\cdot)$ evaluated at the KKT point from Definition~\ref{KKTdef} with respect to their specified argument.  Note that $Z_k$ represents the $k^{\text{th}}$ row of the matrix $Z$.

The proposed prices \eqref{pprice}-\eqref{hprice} are generally discriminatory because they are different for each generator $j$.  In other words, they are not generally locational prices and two generators at the same bus could see different prices.  However, if the input features are uniquely associated with each bus, e.g. $(\bold{p},Z\bold{h})$, then the function $f(\cdot)$ is not a function of $\bold{r}$ or $\bold{g}$, the gradient of $f(\cdot)$ with respect to $\bold{r}_j$ and $\bold{g}_j$ evaluates to $\nabla_{\bold{g}_j} f(\cdot)=\nabla_{\bold{r}_j} f(\cdot)=0$, and the prices become locational prices, or uniform prices.

Discriminatory prices are not well accepted in practical electricity markets for multiple reasons.  Perhaps most importantly, a discriminatory pricing structure lacks transparency and provides poor investment signals as compared to a uniform pricing structure.  For example, ISOs typically provide electricity prices as public information to help guide investment decisions in generation capacity.  In the context of discriminatory pricing, each publicly posted price would be specific to a generator; however, characteristics of each generator would not be provided to the public, due to privacy concerns.  As a result it would not be clear to a potential market entrant why a specific generator is assigned a specific price.  This makes it difficult for a potential market entrant to estimate future profits and justify investment.  In contrast, uniform prices are assigned to a specific location and a small market entrant can accurately estimate future profits using a price-taker assumption.

\subsection{Incentive Alignment}\label{sec:incentives}
Incentive alignment occurs when all market participants have incentive to behave in a way that minimizes total system costs.  In our setting, this occurs when generators have incentive to follow their dispatch instructions, represented by $(\bold{g}_j^{\star},\bold{r}_j^{\star})$.  The following proposition states that the prices \eqref{pprice}-\eqref{hprice} support the dispatch by ensuring that the generator dispatch solves the generator maximization problems \eqref{ProfMax}.

\Prp{\label{Prop1}Consider a KKT point from Definition~\ref{KKTdef}. If prices are chosen as in \eqref{pprice}-\eqref{hprice} and reserve prices are non-negative $\alpha\geq \bold{0}$, then the dispatch $(\bold{g}_j^{\star},\bold{r}_j^{\star},\bold{h}_j^{\star})$ represents a globally optimal solution to the profit maximization~problems~\eqref{ProfMax}. \label{thm1}} \normalfont

Sketch of \proofname:  
The KKT conditions of the TSC-OPF problem \eqref{TSCOPF} imply the KKT conditions of the profit maximization problems \eqref{ProfMax}.  If the reserve prices are non-negative, then the KKT conditions of \eqref{ProfMax} imply the KKT conditions of a convex relaxed version of \eqref{ProfMax} where \eqref{nonsmoothconstraint} is replaced by the inequality constraint \mbox{$\bold{h}_j \leq \min(\bold{\bar{h}}_j,\bold{\bar{g}}_j-\bold{g}_j)$}.  Since the dispatch solves the KKT conditions for this convex relaxation of \eqref{ProfMax}, it must be globally optimal for this convex problem.  Since the dispatch is also feasible for the original non-convex version of \eqref{ProfMax}, it must be globally optimal for \eqref{ProfMax}. \qed

Due to the non-smooth reserve constraint \eqref{nonsmoothconstraint}, this result requires the assumption that the reserve price is non-negative.  Of coarse, this assumption holds when reserve is not used as an input to the NN surrogate, in which case $\nabla_{\!Z_k \bold{h}} f^{\star}=0$.  This assumption is also satisfied for all TSC-OPF instances analyzed in our numerical results.  However, this assumption may not always be satisfied because the value of $\nabla_{\!Z_k \bold{h}} f^{\star}$ may potentially be negative.  Future work will focus on removing this assumption by relaxing the non-smooth constraint.

Proposition \ref{Prop1} requires the dispatch to represent a KKT point as in Definition \ref{KKTdef}.  This does not necessarily require that the dispatch is globally optimal.  In fact, the KKT point may represent a local minimizer, local maximizer, or a saddle-point.  In any case, this incentive alignment result still holds.   

\subsection{Choice of NN Input Features}

In general, the NN surrogate model is likely to be more accurate and better represent stability if the input features contain a full characterization of the system state.  In this context, it is important to note that inputs $(\bold{g},\bold{d},\bold{r},\bold{l})$ are sufficient to uniquely determine the remaining inputs $(Z\bold{h},\bold{p},\bold{q})$, given the fixed parameters of the TSC-OPF problem.  On this basis, it could be expected that the set of inputs $(\bold{g},\bold{d},\bold{r},\bold{l})$ yield the most accurate NN surrogate model.

On the other hand, when the TSC-OPF problem is used to set prices and clear an electricity market, restricting the choice of inputs to those that yield uniform pricing structures is often more desirable.  This observation leads to a trade-off between accuracy and sensible pricing structures that is further analyzed in the numerical results.

\section{Active Sampling Algorithm}\label{SecIII}
Developing an accurate NN model requires training data to be thoroughly sampled along the stability boundary of feasible dispatch values.  Obtaining these samples by sampling the feasible region of the AC-OPF problem, without considering the underlying dynamic models, requires a prohibitively large number of samples to generate and train with.  Rather, this section contributes a model-based active sampling procedure that generates samples along the stability boundary. To identify such samples an Active Search (AS) OPF problem is formulated and solved during each iteration of the algorithm.  Furthermore, this active sampling approach incorporates random sampling of fixed load to ensure the resulting NN is accurate for varying TSC-OPF problem instances.

For the illustration purposes of this paper, a stability constraint is trained for a contingency denoted by $\mathcal{C}$.  This contingency can represent a generator outage, line outage, three phase fault, etc.  Furthermore, the NN $f(\cdot)$ is defined to represent a binary classifier.  In particular, the NN $f(\cdot)$ is assumed to output a scalar that falls in the interval $[0,1]$.  This NN is trained using binary output samples where a classification of $1$ represents a stable solution and a classification of $0$ represents an unstable solution. In this context, the constant $c\in[0,1]$ in constraint \eqref{stabilityconstraint} represents a \emph{threshold parameter} that increasingly encourages stability outcomes as it approaches the value of one.  Finally, in this section, the definition of stability is general but it is inherently assumed that stability is evaluated using dynamic simulation. 

\subsection{Algorithm Intuition}
The proposed model-based active sampling algorithm generates samples in feasible regions where the stability outcome is uncertain according to the NN surrogate model. Each iteration of the algorithm produces a batch of samples that represent optimal solutions of the AS-OPF problem, which has multiple optimal solutions, and updates the NN surrogate model.  As the NN is updated it becomes more accurate and the generated samples cluster near the stability boundary. The key intuition of the approach is that stability is easy to approximate in a large region of the AC-OPF feasible set and the difficulty lies in approximating the stability boundary.  Thus, the proposed approach uses the \emph{optimization to guide the sampling to those regions of the AC-OPF's feasible space that represent the transition between stability and instability.}

\begin{figure}[t]
\subfloat[Early Iteration]{\label{fig:EarlyIterate}\includegraphics[scale=.25]{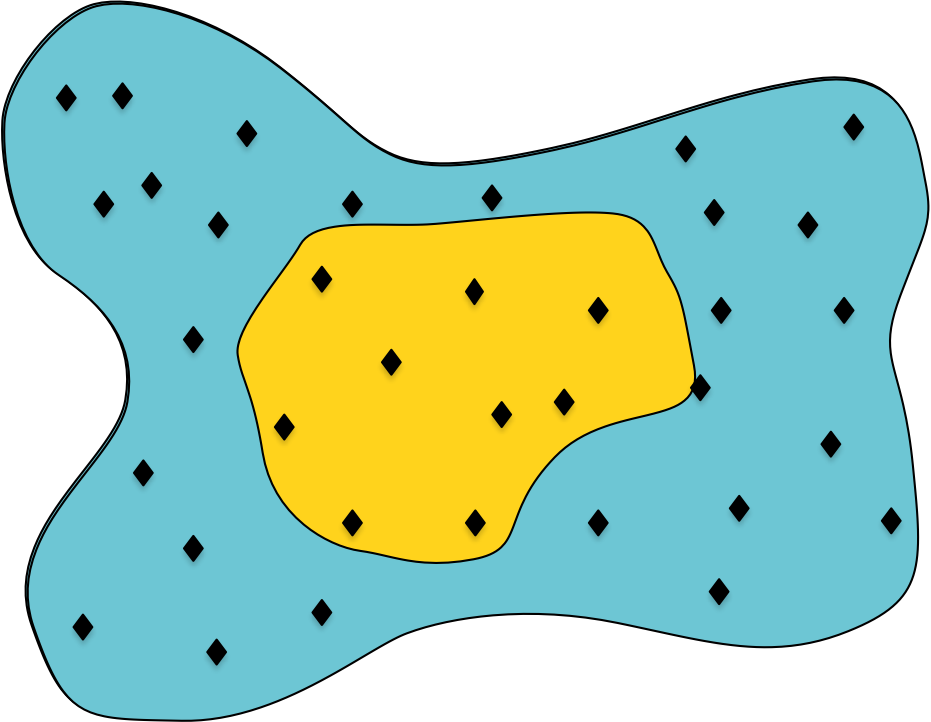}}
\subfloat[ Late Iteration]{\label{fig:LateIterate} \includegraphics[scale=.25]{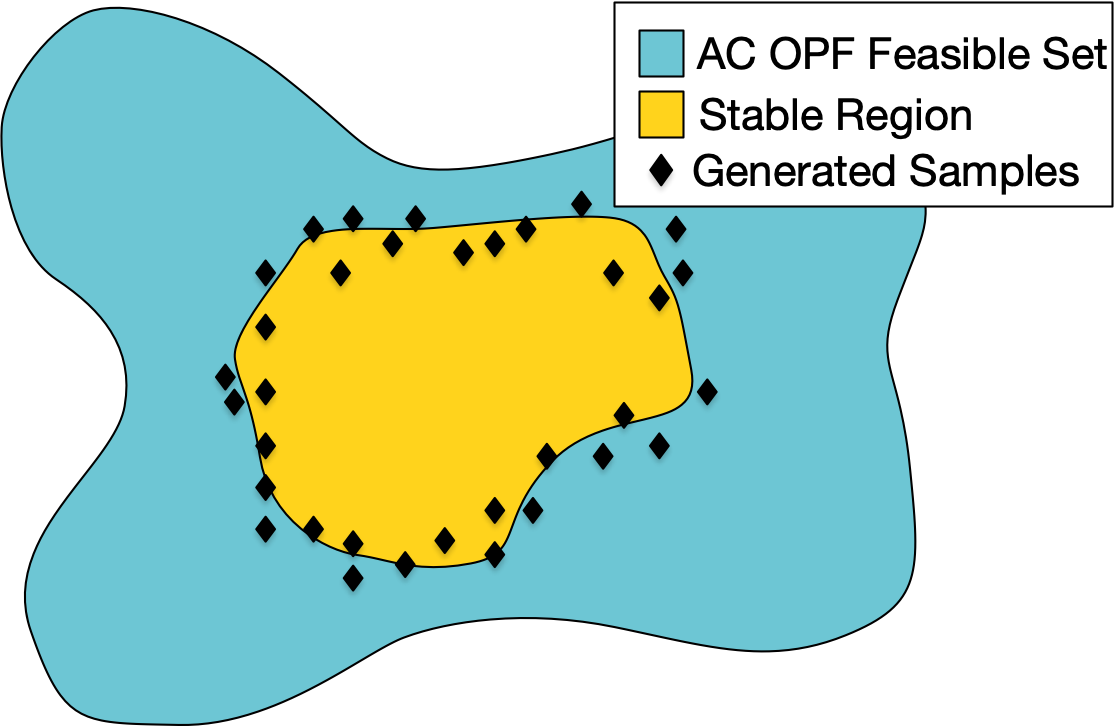}}
\vspace{-5pt}
\caption{\label{fig:ToyExample} Toy example illustrating the intuition of the active sampling algorithm.}
\end{figure}

Figure~\ref{fig:ToyExample} provides a toy example intended to build intuition.  In this toy example the load is considered deterministic and thus load samples need not be generated during each iteration of the active sampling algorithm.  The blue region in Figure~\ref{fig:ToyExample} represents the AC-OPF feasible set and the yellow region represents the set of stable solutions to the AC-OPF problem.  Figure~\ref{fig:EarlyIterate} represents the batch of samples generated by an early iteration of the active sampling algorithm where the NN is not very accurate and most regions in the AC-OPF feasible set have uncertain stability outcomes.  In this case the samples are evenly distributed throughout the AC-OPF feasible set and the algorithm is broadly searching for unstable regions.  Figure \ref{fig:LateIterate} represents the batch of samples generated by a late iteration of the active sampling algorithm where the NN is fairly accurate and the stability boundary is fairly well known.  In this late iteration the generated samples lie very close to the boundary between the stable and unstable regions of the AC-OPF feasible set. 

\subsection{Algorithm Description}
The active sampling algorithm is outlined in Algorithm \ref{ActiveSamplingAlg}.  This algorithm iteratively generates samples that fall near stability boundaries according to the NN trained in the previous iteration.  These new training samples are then added to the existing set of training samples, denoted by $\mathcal{S}$, and the NN is re-trained on the updated set. This iterative process is capable of identifying training samples that fall along the stability boundary and that fall in under-sampled regions of the input space.  

During each iteration of the active sampling algorithm a batch of $\bar{s}$ samples are generated.  These samples are generated in parallel and this generation process is represented by the inner loop beginning at step 3 in Algorithm \ref{ActiveSamplingAlg}.  Each of these samples is generated using the following process:
\begin{enumerate}
    \item Randomly sample the load $(\bold{d},\bold{l})\sim\mathcal{F}$, where $\mathcal{F}$ represents a probability distribution of real and reactive load. This load distribution may represent a probabilistic forecast of the future load or represent a uniform sampling of all possible feasible load realizations.
    \item Solve the AS-OPF problem \eqref{ActiveSearch} to determine input samples. These input samples will represent points in the NN input space that that are near the stability boundary and are AC feasible. (Note: the first iteration $k=1$ solves the AC-OPF problem instead of the AS-OPF problem.)
    \item Run a dynamic simulation of contingency $\mathcal{C}$ that is initialized by the identified solution of the AS-OPF problem and post-process the results to determine output samples.
\end{enumerate}

\begin{algorithm}[t]
\caption{Active Sampling Algorithm} \label{ActiveSamplingAlg}
\begin{algorithmic}[1]
\Require{Samples generated per iteration $\bar{s}$, iterations $\bar{k}$, load distribution $\mathcal{F}$, considered contingency $\mathcal{C}$.}
\State $\mathcal{S}=\emptyset$
\For{$k=1,2,\ldots,\bar{k}$}
\For{$s=1,2,\ldots,\bar{s}$} (parallel computation)\label{line:InputSamples}
\State Randomly sample load $(\bold{d},\bold{l})\sim\mathcal{F}$
\If{k=1} (initial sampling with no NN)
    \State Solve AC-OPF with $(\bold{d},\bold{l})$, e.g. \eqref{TSCOPF} w/o \eqref{stabilityconstraint}   
    \State Collect input sample $I^{(s)}\!\subset\!(\bold{g}^{\star}\!,\!\bold{d},\!\bold{r}^{\star}\!,\!\bold{l},Z\bold{h}^{\star}\!,\!\bold{p}^{\star}\!)$
\Else
    \State Solve AS-OPF \eqref{ActiveSearch} with $(\bold{d},\bold{l})$ and $f^{(k-1)}$
    \State Collect input sample $I^{(s)}\!\subset\!(\hat{\bold{g}},\bold{d},\hat{\bold{r}}\!,\bold{l},Z\hat{\bold{h}},\hat{\bold{p}})$
\EndIf
\State Initialize state of dynamic simulation using $I^{(s)}$
\State Run dynamic simulation of contingency $\mathcal{C}$
\State Collect output sample $O^{(s)}$ (e.g. stability indicator)
\State Store training sample $\mathcal{S}=\mathcal{S}\cup \{(I^{(s)},O^{(s)})\}$
\EndFor
\State Train NN $f^{(k)}$ on current set of training samples $\mathcal{S}$
\label{line:initialization_decay_mu}
\EndFor
\end{algorithmic}
\end{algorithm}

\subsection{Active Search OPF}\label{sec:ASOPF}
Central to the active sampling algorithm is the Active Search OPF problem that determines the input samples in \mbox{steps 9-10} of Algorithm \ref{ActiveSamplingAlg}.  This optimization problem is identical to the AC-OPF problem but with an objective function that minimizes $(f(\cdot)-0.5)^2$.  Intuitively, this objective function is designed to maximize output uncertainty of the NN using a smooth function that is well suited for interior-point algorithms.  

In the active sampling algorithm, the solution of the AS-OPF problem results in a dispatch that is used to initialize a dynamic simulation.  This dispatch must satisfy the transmission constraints used in the dynamic simulation.  If not, then the dynamic simulation will not be initialized in steady state and will instead experience transients before the contingency has even occurred.  For this reason, the AS-OPF problem must optimize over the same transmission model used in the dynamic simulation, which typically represents the AC power flow physics.  This is why this problem is posed as a variant of the AC-OPF problem.  

The AS-OPF problem is written as follows:
\begin{subequations}\label{ActiveSearch}
\begin{align} 
& \min \  \left(f(\bold{g},\bold{d},\bold{r},\bold{l},Z\bold{h},\bold{p})-0.5\right)^2 \tag{\ref{ActiveSearch}}\label{TSCOPFactive}\\
&  \hspace{20pt} \bold{d}-M\bold{g}+\bold{p}= \mathbf{0}
\label{rpowerbalanceactive}\\
& \hspace{20pt} \bold{l}-M\bold{r}+\bold{q}= \mathbf{0} 
\label{qpowerbalanceactive}\\
&\hspace{20pt}(\!\bold{g}_{j}\!,\bold{r}_{j}\!,\bold{h}_{j}\!)\in\mathcal{X}_j \hspace{5pt} \forall j \in \mathcal{G}\\
&\hspace{20pt}(\bold{p},\bold{q})\in\mathcal{T}
\end{align}\end{subequations}
where a general solution to this problem is denoted by $(\hat{\bold{g}},\hat{\bold{r}}\!,\hat{\bold{h}},\hat{\bold{p}},\hat{\bold{q}})$.

There are often many optimal solutions to this problem, each attaining an objective of zero.  For this reason, it is important for the active sampling algorithm to insert randomness into the optimization solver to encourage diverse samples.  This work uses an interior-point algorithm initialized at a randomly sampled starting point to address this problem.

\section{Experimental Setting}
This section describes the experimental setting for the numerical results. A specific frequency stability metric is used to demonstrate the approach.  In particular, the dynamic power system is considered stable in response to a contingency if the frequency at each bus remains above a specified critical frequency $\omega_{\min}$.  Additionally, a specific contingency is considered, namely, the outage of the largest generator in the system.  This stability metric and contingency are chosen because they have been used to construct stability constraints by the ISO in Texas~\cite{li2018design,liu2018participation}. Future work will extend this analysis to multiple contingencies and stability metrics.

\subsection{Turbine Governor Limits}\label{Sec:GovLimits}
Some ISOs today enforce stability constraints in the form of reserve constraints where reserve is procured for ancillary services that contribute to primary frequency control~\cite{mago2020new,li2018design,liu2018participation}.  In this context, reserve ensures each generator $j$ has sufficient headroom $\bar{\bold{g}}_j-\bold{g}_j$ to deploy their automatic transient response.  Procuring reserve in this manner inherently assumes a generator's transient response does not exceed its maximum output limit $\bar{\bold{g}}_j$ during system transients.  This is consistent with the droop requirements for some ISOs.  For example, the ISO in Texas does not require a generator to provide droop control when operating close to its capacity (or High Sustainability Limit)~\cite{ERCOTGuide}.  In contrast, many standard dynamic test cases assume that generators can provide droop control well beyond their capacity, nearly removing capacity limits during system transients.  As a result, in most standard test cases, frequency dynamics do not heavily depend on initial generator output levels, as acknowledged in~\cite{zhang2021encoding}. In an attempt to better represent these maximum output limitations during system transients, turbine governor output limits are set equal to the generators' capacity $\bar{\bold{g}}_j$.  For this reason, the frequency dynamics significantly depend on generator output levels in the numerical results in this paper.

\subsection{Hawaii Test Case}\label{sec:testcase}
A small Hawaii test case~\cite{TexasAMRepo} is used that has 37 buses, 44 generators, and 1136 MW of load.  The test case is altered in two important ways.  First, the three largest generators are removed from the test case.  The considered contingency represents the outage of the largest remaining generator and is referred to as generator 41.  This generator has a maximum output of 87.2 MW and a minimum output of 20 MW.  This generator is also relatively inexpensive and is typically dispatched to its maximum output.

The second important change pertains to generator turbine governor limits. In this test case each generator has a turbine governor of type ``IEEEG1'' or ``GGOV1''~\cite{PowerWorld}. The ``GGOV1'' model always stays within the generator capacity during system transients; however, the ``IEEEG1'' governor models often greatly exceed the generator capacities.  The ``IEEEG1'' model instances are adjusted by setting the maximum turbine governor outputs to represent the generators capacity $\bar{\bold{g}}_j$.

The minimum frequency is defined to be $\omega_{\text{min}}=58.5$ Hz at each bus in the system.  Using all settings specified, for some load profiles this test case produces some AC-OPF solutions that are unstable and some that are stable. Furthermore, the reserve $\bold{h}_j$ is assumed to accommodate droop control for each traditional generator $j$.  In this context, reference~\cite{garcia2021requirements} derives the following expression for the maximum reserve a generator $j$ is capable of providing: 
$$\bold{\bar{h}}_j=\frac{\bar{\bold{g}}_j(\omega_0-\omega_{\min})}{\bold{\nu}_j\omega_0}$$
where $\bold{\nu}_j$ is the droop constant for generator $j$, $\omega_0=60$ Hz is the nominal frequency, and $\omega_{\min}=58.5$ Hz is the minimum frequency.  Intuitively, this value represents the droop control reference signal evaluated at the maximum possible frequency deviation.  It is only an approximation of the amount of droop reserve a generator is capable of providing to arrest frequency decline following a large generator outage. 

\subsection{Choice of Input Features}\label{sec:inputchoice}
Table \ref{InputTable} provides the sets of inputs considered in the remainder of this paper. Input A includes the complex power generation/consumption of each generator/load in the system.  Input A is expected to provide nearly complete information required to determine stability.  Input B then removes reactive power information as compared to Input A.  Input B is expected to result in similar approximation accuracy as Input A because reactive power is weakly coupled with frequency dynamics and the NN is easier to train due to a lower dimensional input space.  Input C further reduces the set of inputs to net real-power injections at each bus, denoted $(\bold{p})$.  Indeed, these input features are expected to result in apparently worse approximation accuracy.  Finally, Input D adds reserve as an input feature, which increases the amount of information available to the NN, and should thus make it more accurate. 

\setlength{\tabcolsep}{3pt} 
\renewcommand{\arraystretch}{1.3} 
\begin{table}[h!]
\vspace{-2pt}
\caption[]{\label{InputTable} Table of proposed inputs for function $f(\cdot)$.}
\vspace{-15pt}
\begin{center}
\begin{tabular}{|c|c|c|c|c|} 
	\hline
	Inputs A&Inputs B&Inputs C&Inputs D\\
	\hline
    $(\bold{g},\bold{d},\bold{r},\bold{l})$&$(\bold{g},\bold{d})$&$(\bold{p})$&$(\bold{p},\bold{1}^T\bold{h})$\\
	\hline
\end{tabular}
\end{center}
\vspace{0pt}
\end{table}

Consistent with the discussion in Section \ref{Sec:Pricing}, some of these choices of inputs result in discriminatory pricing structures (a price for each generator) and some result in uniform pricing structures (a locational price for each bus).  In particular, Inputs A and B result in discriminatory pricing, which is not well accepted in current electricity market practices.  In contrast, Inputs C and D result in uniform pricing structures.  Furthermore, Input D effectively introduces a reserve product into the electricity market and there is only one reserve price seen by all generators in the system.

\subsection{Active Sampling Algorithm Implementation}\label{sec:ASdetails}
The active sampling algorithm was implemented on a Windows-based virtual machine with 128 GB of RAM and 32 2.9 GHz processors each containing two cores.  The inner loop of the active sampling algorithm was parallelized using the PowerShell command ``ForEach-Object'' using 50 concurrent threads. 

The active sampling algorithm was used to train a NN for each of the four choices of inputs separately.  The NN structure was the same for each choice of inputs and included fully connected hidden layers with 128 neurons.  This resulted in a number of parameters that varied between 20,000 and 35,000 for Inputs A - D.  The first and second hidden layers utilized hyperbolic tangent and softplus activation functions, respectively.  The output layer then consisted of a single neuron and used the sigmoid activation function, which scales the NN output to fall between zero and one.  For each choice of inputs, $\bar{s}=500$ samples were collected per active sampling iteration and $\bar{k}=30$ iterations were executed resulting in 15,000 total training samples.  Table \ref{ASTable} provides the time required to perform the active sampling procedure for each input.  Input~D required more time because the AS-OPF problem is more difficult to solve due to the non-smooth constraints \eqref{nonsmoothconstraint} that are enforced.  
\setlength{\tabcolsep}{3pt} 
\renewcommand{\arraystretch}{1.3} 
\begin{table}[h!]
\vspace{-2pt}
\caption[]{\label{ASTable} Time to execute 30 active sampling iterations in minutes.}
\vspace{-15pt}
\begin{center}
\begin{tabular}{|c|c|c|c|c|c|} 
	\hline
	&Inputs A&Inputs B&Inputs C&Inputs D\\
	\hline
	Total Time (min.)&576&555&556&810\\
	\hline
	Mean Time per iteration (min.)&19.2&18.5&18.5&27.0\\
	\hline
\end{tabular}
\end{center}
\end{table}

\subsubsection{Load Sampling}\label{Sec:LoadSampling}
\begin{figure}[h!]
    \centering    \includegraphics[width=0.85\columnwidth]{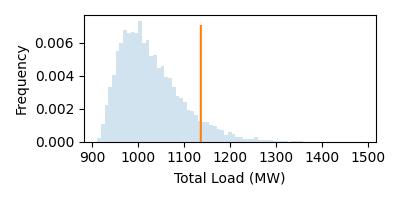}
    \caption{Histogram illustrating the distribution of total system load drawn from a Gamma distribution.}\label{loadhist}
\end{figure}

The load in units of MW was randomly sampled using two steps.  The first step randomly samples the total system load from a Gamma distribution with shape parameter of 3, a shift parameter of 909 MW, and a scale parameter of 40.  To achieve this level of total system load all load is scaled proportionally to its original value provided by the test case.  Real and reactive power are both scaled by the same scaling factor.  The histogram in Figure \ref{loadhist} illustrates 10,000 random samples of total system real power load drawn from this distribution.  The red vertical line marks the total system load for the original test case.  

The second step in the load sampling process adds a randomly sampled perturbation to the real and reactive power consumption of each individual load.  This allows the load profile to change in space across the transmission network.  The perturbation was randomly sampled from a Gaussian distribution centered at zero and with a standard deviation representing 20\% of the load value after being scaled in the first step.

\subsubsection{Solving Optimization}
The Julia package PowerModels.jl was used to formulate the AC-OPF problem in JuMP.  To obtain the TSC-OPF and AS-OPF problems the Julia package MathOptAI.jl was used to embed our NN surrogate models into the AC-OPF problem~\cite{8442948,MathOptAI2024,Lubin2023}.  The optimization problems were solved using the IPOPT interior-point algorithm with MA27 as the linear solver~\cite{wachter2006implementation}.  The default settings in IPOPT were used for all problem instances except for the instances that included reserve, e.g. Input D.  When using Input~D, constraint \eqref{nonsmoothconstraint} must be enforced, which is non-smooth and causes potential numerical errors and convergence issues in IPOPT.  In this case convergence tolerance levels were slightly reduced in IPOPT and the option to compute approximate hessians was used rather than computing exact hessians.  

\subsubsection{Simulations}
The PowerWorld Simulator software was used to run dynamic simulations along with EasySimAuto, a Python package that interfaces with PowerWorld Simulator \cite{PowerWorld,ESA}. The outage of generator 41 was simulated over a horizon of 300 seconds and trajectories were post-processed to determine the lowest frequency.  To reduce the computational burden, an iterative procedure was implemented that terminates once the frequency drops below 58.5 Hz or once the frequency has begun to rise, reducing the simulation horizon.

\subsubsection{NN Training}\label{Sec:Training}
The NN surrogate models were trained using a Stochastic Gradient Descent (SGD) algorithm implemented with PyTorch~\cite{paszke2019pytorch}.  Each time a NN was trained, the set of total samples was partitioned into 80\% training and 20\% validation samples used to define the stopping criteria.  Additionally batch sizes of 64 samples were used.

\section{Numerical Results}\label{SecIV}
This section presents numerical results that validate the performance of the TSC-OPF problem over 1000 random load samples by evaluating stability of the optimal dispatch.  Additionally, trends in the optimal solution, including the number of stable solutions, are observed as the threshold parameter $c$ varies.

With the exception of Sections \ref{Sec:cpu} and \ref{Sec:Feasibility}, all remaining subsections analyze optimal solutions of the TSC-OPF problem, inherently assuming a solution exists.  For this reason the analysis focuses on a subset of threshold values for which a significant number of the 1000 load samples result in a TSC-OPF problem that is found to have a feasible point by IPOPT. In total, there were 960 out of 1000 TSC-OPF instances that were locally solved and feasible for all optimization solves.  Of these 960 load samples, the AC-OPF solution results in 123 unstable dynamic simulations, or 12.8\% of the samples.

\subsection{Validation: Ensuring Stability with Different Input Features}
\begin{figure}[h!]
    \centering
    \includegraphics[width=0.85\columnwidth]{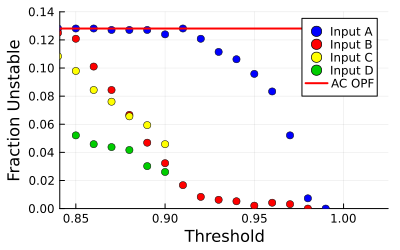}
    \caption{Fraction of unstable load samples versus threshold parameter $c$ for different input selections.}
    \label{fig:stabilityvthreshold1}
\end{figure}

Figure \ref{fig:stabilityvthreshold1} shows the fraction of unstable load samples versus the chosen threshold parameter for different choices of input features.  The red line represents the number of unstable load samples resulting from the solution of the AC-OPF problem, with a value near $0.128 (12.8\%)$.  All choices of inputs result in a fraction of unstable load samples that decreases with the threshold parameter $c$. Inputs A and B achieve stability for all load samples at a threshold parameter of $0.99$ and $0.98$, respectively. Furthermore, Input B achieves fewer unstable load samples for a wider range of threshold values, making it the best selection of inputs.  Inputs C and D reduce the number of unstable solutions to $4.6\%$ and $2.6\%$ respectively, which occurs at a threshold of $c=0.9$.  No data points are provided for Inputs C or D beyond $c=0.9$ because threshold values higher than this cause a significant number of infeasible TSC-OPF problems.

\subsection{TSC-OPF with Input B}
This subsection analyzes trends of solutions to the TSC-OPF problem when using Input B.  Input B is chosen for this analysis because the previous subsection found that it was the best at reducing instability.  

\subsubsection{Impacts on the Optimal Objective Value}
\begin{figure}[h!]
    \centering
    \includegraphics[width=0.85\columnwidth]{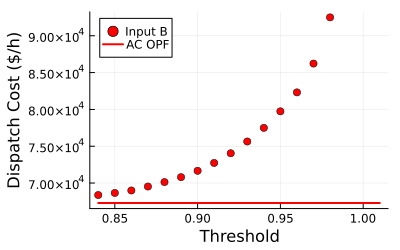}
    \caption{Optimal objective value of TSC-OPF (dispatch cost) versus threshold parameter $c$.}
    \label{fig:ObjvThreshold}
\end{figure}

Increasing the threshold parameter, $c$, restricts the feasible set of the TSC-OPF problem \eqref{TSCOPF} and increases the optimal objective value, representing total dispatch cost. Figure \ref{fig:ObjvThreshold} shows the average total cost versus the chosen threshold parameter.  The red line represents the average dispatch cost produced by the AC-OPF problem.  Costs increase super-linearly with the threshold parameter $c$.

\begin{figure}[h!]
    \centering
    \includegraphics[width=0.85\columnwidth]{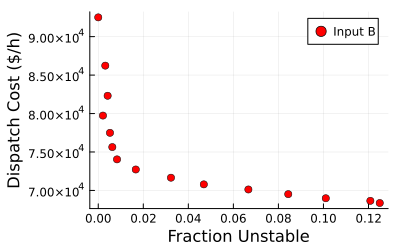}
    \caption{
    \label{fig:ObjvStability1}Optimal objective value of TSC-OPF (dispatch cost) versus the fraction of unstable load samples.}
\end{figure}

When comparing choices for $c$, this threshold is interpreted as a parameter that encourages stability as it approaches $c=1$.  To analyze the trade-off between stability and cost, Figure \ref{fig:ObjvStability1} plots the average objective value versus the fraction of unstable load samples as the threshold parameter varies. For example, a fraction of unstable samples of 0.0 can be obtained, nearly guaranteeing stability, at the expense of increasing cost by approximately 37\% as compared to the cost obtained by the AC-OPF problem.  This performance is achieved by setting the threshold value to $c=0.98$.  Alternatively, the threshold parameter may be chosen to attain a value near the elbow of this curve.  For example, a threshold of $c=0.92$ results in a fraction of unstable load samples of $0.008$ and only increases the optimal objective value by 10\% as compared to the optimal value of the AC-OPF problem. 

\subsubsection{CPU Run Times}\label{Sec:cpu}
The computational effort required to solve an individual TSC-OPF problem was minimal.  Table \ref{Table:CompTimes} provides statistics on computational time for the TSC-OPF problem with Input B and for the original AC-OPF problem.  For each of the 2000 optimization problems, the termination status of IPOPT was either ``Locally Solved'' or ``Locally Infeasible.''  For the TSC-OPF problem a threshold of 0.98 is used, which was chosen to be restrictive but not too restrictive such that many problem instances become infeasible.  The same 1000 validation load samples are used for both problems.  The AC-OPF problem was very fast to solve and solved faster than 1 second on average.  Input B was slower to solve, but only increased the solve time by one order of magnitude.  
\setlength{\tabcolsep}{1pt} 
\renewcommand{\arraystretch}{1.3} 
\begin{table}[h!]
\vspace{-2pt}
\caption[]{\label{Table:CompTimes} Table of Computational Times (1000 Load Samples)}
\vspace{-15pt}
\begin{center}
\begin{tabular}{|c|c|c|c|c|c|c|} 
    \hline
    &$c$&\# Infeas.&Mean Time (s)&Max Time (s)&Mean Iters& Max Iters\\
	\hline
    Input B&0.98&2&8.19&25.77&483.66&1091\\
	\hline
    AC-OPF&n/a&0&0.33&4.02&31.89&44\\
    \hline
\end{tabular}
\end{center}
\vspace{-10pt}
\end{table}

\subsubsection{TSC-OPF Feasibility}
\begin{figure}[h!]
    \centering
    \includegraphics[width=0.8\columnwidth]{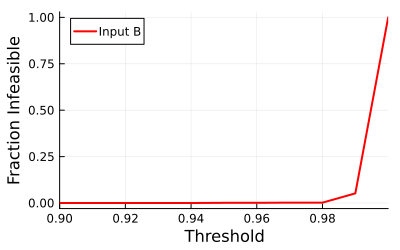}
    \caption{Fraction of infeasible load samples versus threshold parameter $c$.}
    \label{fig:InfvThreshold}
\end{figure}

As the threshold parameter increases, an increasing number of the 1000 load samples become infeasible, as illlustrated in Figure~\ref{fig:InfvThreshold}.  When using Input B, the TSC-OPF problem is capable of accommodating very high threshold values of $c>0.99$ before resulting in a significant number of infeasible TSC-OPF problem instances.  Furthermore, the TSC-OPF problem achieves impressive stability performance for thresholds as low as 0.92, according to Figure \ref{fig:stabilityvthreshold1}. This provides a wide range of threshold values that can be chosen to both ensure stability and feasibility.

\subsubsection{Threshold impacts on frequency trajectory}\label{Sec:Feasibility}
\begin{figure}[h!]
    \centering
    \includegraphics[width=0.85\columnwidth]{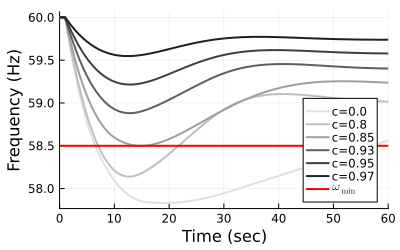}
    \caption{Frequency trajectories with increasing threshold parameter.}
    \label{fig:multfreq}
\end{figure}

This subsection provides an example where the AC-OPF problem results in an unstable solution and the TSC-OPF problem with Input B makes the system stable.  Furthermore, the frequency trajectory is observed as the threshold parameter $c$ varies.  For a specific load sample Figure \ref{fig:multfreq} provides frequency trajectories for thresholds varying from 0 to 0.97.  In this figure, each trajectory represents the averaged frequency at each bus in the system.  Although the frequency limit of 58.5 Hz is enforced at each bus in the system, the frequencies at each bus exhibit very little difference in the small Hawaii test case considered.

The frequency trajectory when $c=0$ corresponds to that produced by the solution of the AC-OPF problem.  The lowest point along this frequency trajectory, also referred to as the \emph{frequency nadir}, is far below 58.5 Hz.  The frequency nadir increases with the threshold parameter, reaching 58.5 Hz when the threshold is approximately $c=0.85$.  As the threshold continues to increase, the frequency nadir  continues to increase to conservatively high values.  

\subsubsection{Impacts on Generator Dispatch}
\begin{figure}[h!]
    \centering
    \includegraphics[width=0.75\columnwidth]{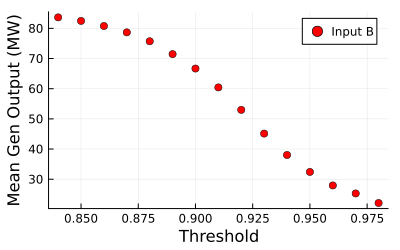}
    \caption{Pre-contingency output of outed generator.}
    \label{GenvThreshold}
\end{figure}

The most effective way to reduce the impact of the generator outage is to reduce the pre-contingency real-power output of the outed generator. Figure \ref{GenvThreshold} illustrates how the TSC-OPF problem dispatches generator 41 downward as the threshold increases.  When the threshold is low, generator 41 is dispatched at its maximum output of 87.2 MW and when the threshold is high it is dispatched at its minimum output of 20 MW. 

\subsection{Comparison to Simple Sampling}
To illustrate the need for the active sampling algorithm, it is compared to a \emph{simple sampling} procedure.  Specifically, using Input B, a NN was trained using 100,000 training samples each generated by randomly sampling the load as in Section \ref{Sec:LoadSampling} and then solving an AC-OPF problem to determine the generation.  These training samples were collected using the same computing platform as stated in Section \ref{sec:ASdetails} and computation was executed on 50 parallel threads.  It required approximately 500 minutes to collect 100,000 samples, which is comparable to the times taken to perform the active sampling algorithm, shown in Table \ref{ASTable}.  After collecting the training samples, the NN was trained using the same procedure as in Section \ref{Sec:Training}.

\begin{figure}[h!]
    \centering
    \includegraphics[width=0.85\columnwidth]{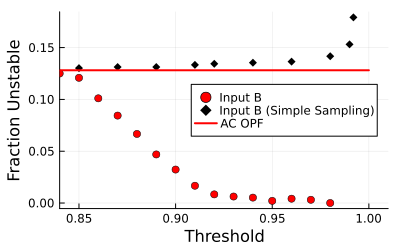}
    \caption{Stability versus threshold parameter $c$ using different sampling procedures.}
    \label{fig:StabilityvThreshold2}
\end{figure}

Figure \ref{fig:StabilityvThreshold2} illustrates that the TSC-OPF problem is unable to improve stability when using the NN surrogate model trained by the simple sampling procedure.  In this case the TSC-OPF problem begins to cause more instability at high threshold values indicating that it is poorly approximating the stability boundary.  This demonstrates the contribution of our active sampling algorithm.

\subsection{TSC-OPF with Inputs C and D}
Figure \ref{fig:stabilityvthreshold1} shows that Inputs C and D do not perform as well as Input B in terms of stabilizing the system.  However, Section \ref{Sec:Pricing} points out that Input B results in a discriminatory pricing structure that is not well accepted in an electricity market setting.  As a result, an ISO may prefer to implement Inputs C or D because they result in a uniform pricing structure that satisfies an incentive alignment property.  

\begin{figure}[h!]
    \centering
    \includegraphics[width=0.85\columnwidth]{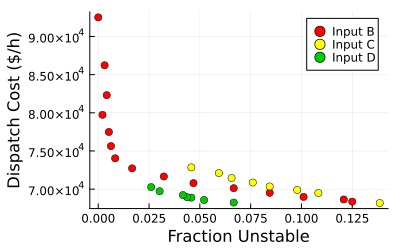}
    \caption{Optimal objective value of TSC-OPF (dispatch cost) versus the fraction of unstable load samples.}
    \label{fig:ObjStability}
\end{figure}

Figure \ref{fig:ObjStability} compares the performance of Inputs B, C, and D in terms of both stability and cost.  Input D outperforms Input C in terms of both cost and stability.  Input D is not able to attain low instability numbers like Input B; however, it is able to achieve low cost solutions at reasonably low instability numbers. In fact, Input D is able to achieve a fraction of stability of $0.026$ when increasing the cost by only 7\%.

The prices are not analyzed when using Input B because it is not clear how to set prices in a way that is well accepted in practice.  Instead, prices are compared when using Inputs C and D.  Both Inputs C and D result in locational prices for real and reactive power.  Additionally, only Input D has a price for reserve, because the reserve price for Input C is effectively zero.  

\begin{figure}[h!]
    \centering
    \includegraphics[width=0.75\columnwidth]{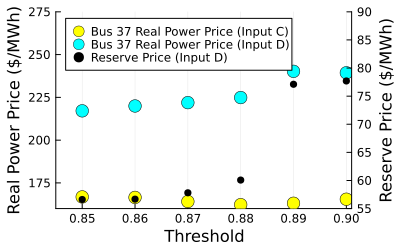}
    \caption{Average prices seen by outed generator during extreme load samples.}
    \label{prices}
\end{figure}

Figure \ref{prices} illustrates pricing trends at bus 37 during highly stressed load samples.  This figure focuses on the 7\% of load samples that result in the highest dual variable $\gamma$, using this as a metric for load samples that are the most difficult to stabilize.  Using these load samples, Figure \ref{prices} plots the average real-power price at bus 37 and the average reserve price for different choices of the threshold parameter.  Notice that the real-power price at bus 37 is much lower when using Input C as compared to Input D. This is because the TSC-OPF problem must dispatch generator 41 downward to maintain stability and must choose prices/dispatch values in a way that allows generator 41 to maximize its profit, see Proposition \ref{Prop1}.  For this reason the real power price at bus 37 must remain lower.  This will also impact other generators at bus 37, effectively reducing their profits.

Unlike Input C, Input D is capable of dispatching generator 41 downward by increasing the price of reserve. Furthermore, reserve prices tend to increase with the threshold parameter, encouraging generator 41 to increase its reserve $\bold{r}_j$ and reduce its generation output $\bold{g}_j$.  This alternative method of dispatching generator 41 downward allows Input D to retain higher real-power prices at bus 37 during stressful load samples.  As a result, when using Input D other generators at bus 37 will continue to see high profits during these stressful load samples.  Interestingly, generator 41 will also see higher profits when using Input D.  This is because prices remain high at bus 37 and because generator 41 will receive payment for providing reserve even though it is not available after being outed to provide that reserve.

\section{Conclusion}\label{SecV}
This paper developed a TSC-OPF problem that enforces frequency stability constraints using NN surrogate models.  An active sampling algorithm was developed that iteratively generates NN input samples along the stability boundary while remaining AC feasible.  In the context of wholesale electricity markets, a pricing structure was presented using a general set of NN inputs that were shown to satisfy an incentive alignment property, ensuring generators have incentive to follow their dispatch instructions.  These prices were shown to be discriminatory for certain choices of NN inputs and uniform for others.  The TSC-OPF solutions were validated over randomly sampled load using a Hawaii test case and were shown to significantly enhance frequency stability at low computational cost and low financial cost.  For certain choices of NN inputs and threshold parameters, the TSC-OPF problem was able to attain a stable dispatch for all load samples, including the 12.8\% of load samples for which the AC-OPF problem resulted in an unstable dispatch. Numerical results show that the choices of NN inputs that most effectively enhance stability also result in discriminatory pricing structures.  Furthermore, there exists choices of NN inputs that are reasonably effective at enhancing stability and result in uniform pricing structures.

\section*{Acknowledgment}
The authors would like to thank Dr. Michael Grosskopf for discussions and valuable suggestions.

\bibliography{bibfile}
\bibliographystyle{IEEEtran}

\end{document}